\documentclass[12pt,a4paper]{article}

\usepackage{amsmath,amssymb,latexsym,graphics,epsfig}
\usepackage{hyperref}
\usepackage{color}
\usepackage{amsthm}
\usepackage{graphicx,url}
\usepackage{amsopn}
\usepackage{amssymb}
\usepackage[english]{babel}
\usepackage{soul}
\usepackage{blkarray}

\newtheorem{Theorem}{Theorem}[section]
\newtheorem{Lemma}[Theorem]{Lemma}
\newtheorem{Corollary}[Theorem]{Corollary}
\newtheorem{Proposition}[Theorem]{Proposition}
\newtheorem{Definition}[Theorem]{Definition}
\newtheorem{Example}[Theorem]{Example}
\newtheorem{Conjecture}[Theorem]{Conjecture}

\DeclareMathOperator{\Circ}{circ}

\DeclareMathOperator{\diam}{diam}
\DeclareMathOperator{\dist}{dist}

\DeclareMathOperator{\ecc}{ecc}

\DeclareMathOperator{\rad}{rad}

\def\j{\mbox{\boldmath $j$}}

\def\u{\mbox{\boldmath $u$}}
\def\x{\mbox{\boldmath $x$}}
\def\y{\mbox{\boldmath $y$}}
\def\z{\mbox{\boldmath $z$}}

\def\vecrho{\mbox{\boldmath $\rho$}}

\def\u{\mbox{\boldmath $u$}}

\def\v{\mbox{\boldmath $v$}}

\def\vec0{\mbox{\boldmath $0$}}

\def\A{\mbox{\boldmath $A$}}

\def\D{\mbox{\boldmath $D$}}

\begin{document}

\title{On some metric properties \\
of supertoken graphs\thanks{The second and third authors' research has been supported by AGAUR from the Catalan Government under project 2021SGR00434 and MICINN from the Spanish Government under project PID2020-115442RB-I00.
The research of the third author was also supported by a grant from the  Universitat Polit\`ecnica de Catalunya with references AGRUPS-2022 and AGRUPS-2023.}}
			
\author{
	E. T. Baskoro\\
	{\small Faculty of Mathematics and Natural Sciences, Institut Teknologi Bandung} \\
	{\small Bandung, Indonesia} \\
	\vspace{.25cm}
	{\small {\tt{ebaskoro@itb.ac.id}}}\\
C. Dalf\'o\\
		{\small Departament de Matem\`atica, Universitat de Lleida} \\
		{\small Igualada (Barcelona), Catalonia} \\
				\vspace{.25cm}
		{\small {\tt{cristina.dalfo@udl.cat}}}\\
	M. A. Fiol\\
    	{\small Departament de Matem\`atiques, Universitat Polit\`ecnica de Catalunya} \\
    	{\small Barcelona Graduate School} \\
     {\small  Institut de Matem\`atiques de la UPC-BarcelonaTech (IMTech)}\\
    	{\small Barcelona, Catalonia} \\
    	\vspace{.25cm}
    	{\small {\tt{miguel.angel.fiol@upc.edu}}}\\
R. Simanjuntak\\
{\small Faculty of Mathematics and Natural Sciences, Institut Teknologi Bandung} \\
{\small Bandung, Indonesia} \\
\vspace{.25cm}
{\small {\tt{rino@math.itb.ac.id}}}
}
	\date{}

\maketitle

\newpage

\begin{abstract}
In this paper, we construct two infinite families of graphs $G(d,c)$ and $G^+(d,c)$, where, in both cases, a vertex label is 
$x_1x_2\ldots x_c$ with $x_i\in\{1,2,\ldots, d\}$. We provide a lower bound on the metric dimension, tight on $G^+(d,c)$. Moreover, we give the definition and properties of the supertoken graphs, a generalization of the well-known token graphs. Finally, we provide an upper bound on the metric dimension of supertoken graphs. 
\end{abstract}




\section{Introduction}

 Let $G=(V,E)$ be a graph with a (finite) set
$V=V(G)$ of vertices and a set $E=E(G)$ of edges.
The distance from vertex $u$ to vertex $v$ in $G$, denoted by $\dist_G(u,v)$, is the length of a shortest path from $u$ to $v$. So, the diameter $d=\diam(G)$ is the distance between any pair of the furthest vertices in $G$.
If $G$ has $n$ vertices, its distance matrix $\D$ is an $n\times n$ matrix with entries $(\D)_{u,v}=\dist_G(u,v)$.
This matrix has been extensively studied in the literature and has interesting applications, both theoretical and practical (for instance, in telecommunication, chemistry, etc.); see Edelberg, Garey, and Graham \cite{egg76}, and Dededzi \cite{k14}.
The results in the following theorem, which we use in this paper, are due to
Graham and Pollak \cite{gp71}, and  Bapat, Kirkland, and Neumann \cite{bkn05}.
\begin{Theorem}
\label{th:detD}
Let $G$ be a graph with distance matrix $\D(G)$. Then, the following statements hold.
\begin{itemize}
\item[$(i)$]\cite{gp71}
If $G=T$ is a tree, the determinant of $\D(T)$  only depends on its number of vertices $n$: 
$\det \D(T)=(-1)^{n-1}(n-1)2^{n-2}$.
\item[$(ii)$]\cite{bkn05} 
Let $G$ be a unicyclic graph. Then,  $\det \D(G)=0$ when the only cycle of $G$  has an even number of edges; and $\det \D(G)=(-2)^m [k(k + 1) + \frac{2k+1}{2} m]$ if $G$ has $2k+1+m$ vertices and a cycle with $2k+1$ edges. 
\end{itemize}
\end{Theorem}

A vertex subset $C=\{z_1,z_2,\ldots,z_c\}\subset V$ is a {\em resolving set} if 
every vertex $u\in V$ is determined by the vector 
$$
\vecrho=\vecrho(u|C)=(\dist_G(u,z_1),\dist_G(u,z_2),\ldots,\dist_G(u,z_c)),
$$ 
which is called the {\em representation (or position) of $u$ with respect to $C$}.
 The {\em metric dimension} of $G$, denoted by $\dim(G)$, is the minimum cardinality of a resolving set.
 Notice that a vertex subset $C\subset V$ is a resolving set of $G$ if 
the $|C|\times |V|$ submatrix of $\D(G)$, with rows indexed by the vertices of $C$, has all its columns different,
and the {\em metric dimension} of $G$ is the minimum cardinality of such a subset.
The notion of metric dimension was introduced independently by Harary and Melter \cite{hm76} and by Slater \cite{s75}. Since then, there has been extensive literature on the topic.
For example, 
C\'aceres, Hernando, Mora, Pelayo, Puertas, Seara, and Wood \cite{chmppsw07}
studied the metric dimension of the Cartesian product $G\Box H$ of graphs $G$ and $H$. As one of their main results, they provided a family of graphs $G$ with bounded metric dimension for which the metric dimension of $G\Box G$ is unbounded.  Moreover, Bailey and  Cameron \cite{bc11} dealt with the metric dimension of the Johnson and the Kneser graphs. In particular, they computed the metric dimension of the former, which can be defined as the $k$-token graph of $K_n$ (see the end of the following paragraph). 
In general, it is known that computing the metric dimension of a graph is an NP-complete problem; see Garey and Johnson \cite[p. 204]{gj79}.
Some applications of the metric dimension approach are in network theory and combinatorial optimization,
see Chartrand and Zhang \cite{cz03} for a survey.

In this paper, we use a generalization of symmetric powers or tokens graphs, which we call supertoken graphs.  Audenaert, Godsil, Royle, and Rudolph~\cite{agrr07} defined the {\em $k$-symmetric power} of a graph $G=(V,E)$, with vertices the $k$-subsets of $V$, and two vertices are adjacent if their symmetric difference is the end-vertices of an edge in $E$. Fabila-Monroy, Flores-Pe\~{n}aloza, Huemer, Hurtado, Urrutia, and Wood~\cite{ffhhuw12} renamed them as the {\em $k$-token graph} of $G$, and denoted it by $F_k(G)$. The reason is that each vertex of $F_k(G)$ corresponds to a configuration of $k$ tokens placed in the (different) vertices of $G$, and the adjacency is obtained by moving one token from a vertex to an adjacent vertex.  
The Laplacian spectra of token graphs were studied by Dalf\'o, Duque, Fabila-Monroy, Fiol, Huemer, Trujillo-Negrete, and  Zaragoza Mart\'inez \cite{ddffhtz21}.
As an important example of token graphs, we must mention that the Johnson graphs can be defined as $J(n,k)=F_k(K_n)$ for any integers $n\ge 2$ and $1\le k<n$.

This paper is structured as follows. In the following section, we consider the construction and properties of two infinite families of graphs $G(d,c)$ and $G^+(d,c)$, which, in both cases, can be considered as graphs on alphabets. Moreover, we provide a lower bound on the metric dimension, which is tight on $G^+(d,c)$. In Section \ref{sec:supertoken}, there is the definition of the supertoken graphs ${\cal F}_k(G)$, and we find its metric properties when $G$ is a complete graph $K_n$ and a general graph. Finally, in Section \ref{sec:metric-dim}, there is an upper bound on the metric dimension of 
${\cal F}_k(K_n)$, and ${\cal F}_k(G)$ for a general graph $G$.


\section{Maximal graphs with given metric dimension}
\label{sec:resolving}

In this section, we construct two families of graphs on alphabets to obtain graphs with maximum numbers of vertices for a given metric dimension. In general, graphs on alphabets are constructed by labeling the vertices with words on a given alphabet, and
specifying rules that relate pairs of different words to define the edges, see G\'omez, Fiol, and Yebra \cite{gfy92}. 

Specifically, the graphs in the following definition allow us to get a tight lower bound for the metric dimension.
\begin{Definition}
Given integers $d$ and $c$,  the graph $G(d,c)$ has every vertex represented by a sequence $x=x_1x_2\ldots x_c$ with $x_i\in[1,d]=\{1,2,\ldots,d\}$. Moreover, two vertices $x$ and $y=y_1y_2\ldots y_c$ are adjacent if and only if $|x_i-y_i|\in \{0,1\}$ for every $i=1,2,\ldots,c$.
\label{def:G(d,c)}
\end{Definition}
The following result shows some basic metric properties of the graph $G(d,c)$. 
\begin{Lemma}
\label{lem:diam+}
Let $x$ and $y$ be two generic vertices of the graph $G=G(d,c)$. 
\begin{itemize}
\item[$(i)$]
   If $x=x_1x_2\ldots x_c$ and $y=y_1y_2\ldots y_c$, the distance between $x$ and $y$ is
   $$
   \dist_G(x,y)=\max_{i\in[1,c]}\{|x_i-y_i|\}.
   $$
   \item[$(ii)$]
   The eccentricity of the vertex $x$ is 
   $$
  \ecc_G(x)=\max_{i\in[1,c]}\{x_i,d-x_i\}.
   $$
\item[$(iii)$]
   The diameter of the graph $G(d,c)$ is 
   $$
   \diam(G(d,c))=d-1.
   $$
   \end{itemize}
\end{Lemma}
\begin{proof}
$(i)$ From the adjacency conditions in Definition \ref{def:G(d,c)}, it is clear that  $\dist_G(x,y)=\max_{i\in[1,c]}\{|x_i-y_i|\}$. To follow a shortest path from $x$ to $y$, note that each step $x_i$ can be changed to $x_i+1$ or $x_i-1$ for every $i=1,\ldots,c$. Thus, each $x_i$ can be changed to $y_i$ in exactly $|x_i-y_i|$ steps and, hence, $\dist(x,y)$ is as claimed.\\
$(ii)$
By $(i)$, every vertex $y$ satisfies $\dist_G(x,y)\le \max_{i\in[1,c]}\{x_i,d-x-i\}$, and, for example, a vertex at maximum distance from $x$ is either $11\ldots 1$ or $dd\ldots d$.
\\
$(iii)$ This is a consequence of $(i)$ and $(ii)$ since $\max_{i\in[1,c]}\{|x_i-y_i|\}=d-1$ and 
$\min_{i\in[1,c]}\{x_i\}=1$.\\
\end{proof}

\begin{Proposition}
Let $G=(V,E)$ be a graph with $n$ vertices and diameter $d$. Let $c=c(d,n)$ be the minimum integer satisfying $n\le d^c+c$. Then, the metric dimension of $G$ satisfies
\begin{equation}
\label{bound-dim}
\dim(G)\ge c(d,n),
\end{equation}
and the bound is tight.
\end{Proposition}
\begin{proof}
Let  $C\subset V$ be a minimum resolving set with vertices $z_1,z_2,\ldots,z_c$. Since every vertex $u\in V$
must be univocally identified with a sequence $x_1x_2\dots x_c$, where $x_i=\dist_G(u,z_i)\in [0,d]$, it must be
$
n\le d^c+c.
$
If $x_i=0$ for some $i$, then $u=z_i$ and the other values of $x_j$, for $j\neq i$, are determined by $x_j=\dist_G(z_i,z_j)$.
To show that the bound is tight, consider the graph $G^+(d,c)$ 
obtained from $G(d,c)$ by adding the vertices $w_1,w_2,\ldots,w_c$ and, for every $i=1,2,\ldots,c$, the vertices $w_i$ are adjacent to all vertices $x=x_1x_2\ldots x_c$ with $x_i=1$. The graph $G^+(d,c)$ has order $d^c+c$ and diameter $d$, as required. Moreover, from its construction, $\dist_G(x,w_i)=x_i$, so that the set $C=\{w_1,w_2,\ldots,w_d\}$ is a resolving set. 
\end{proof}

As an example, we show the graph $G^+(4,2)$ in Figure \ref{fig:G^+(4,2)}.

\begin{figure}[t]
	\begin{center}
		\includegraphics[width=6cm]{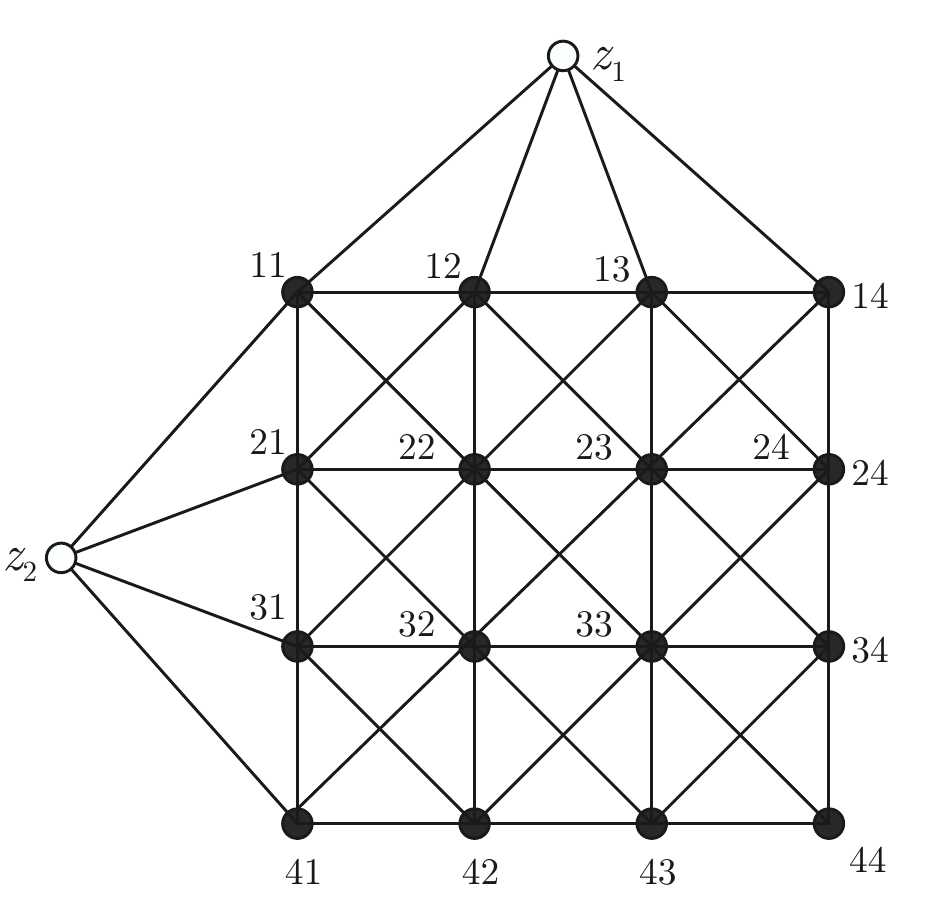}
	\end{center}
 \vskip-.25cm
	\caption{The optimal graph $G^+(4,2)$ with resolving set $\{z_1,z_2\}$, diameter $d=4$ and $d^2+c=18$ vertices.}
	\label{fig:G^+(4,2)}
\end{figure}

For some graphs, the inequality \eqref{bound-dim} can be improved. Let us show an example.
A graph $G=(V,E)$ with diameter $d$ is {\em degree-regular} if every vertex has the same number of vertices at distance $i=1,\ldots,d$. That is, the numbers $k_i(u)=|G_i(u)|$, for $i=0,1,\ldots,d$, do not depend on the vertex $u\in V$, and we simply write $k_i$. Some well-known examples of degree-regular graphs are the vertex-transitive graphs and the distance-regular graphs.
Then, the following result is straightforward. 

\begin{Lemma}
\label{lemma1}
The metric dimension $\mu=\dim(G)$ of a degree-regular graph $G$, with $N$ vertices, diameter $d=d(G)$, and numbers $k_0(=1),k_1,\ldots,k_d$ must satisfy
$$
N\le |{\cal S}(\mu; k_0,k_1,\ldots,k_d)|,
$$
where ${\cal S}(\mu; k_0,k_1,\ldots,k_d)$ is the set of sequences of length $\mu$ having at most $k_i$ numbers $i$, with $i\in [1,d]$.
\end{Lemma}


\section{Supertoken graphs}
\label{sec:supertoken}

Here, we consider the \emph{k-supertoken graph}  ${\cal F}_k(G)$ defined as follows. Let $G$ be a graph on $n$ vertices. Then, for some integer $k\ge 1$, each vertex of ${\cal F}_k(G)$ corresponds to $k$ indistinguishable tokens placed in some of the (not necessarily distinct) $n$ vertices in $G$. Then,  the vertex set of ${\cal F}_k(G)$ corresponds to the combinations with repetitions $CR^n_k$ of $n$ elements taken $k$ at a time, with cardinality
$$
|CR^n_k|= {{n+k-1}\choose{k}} = {{n+k-1}\choose{n-1}}.
$$
Thus, every vertex of ${\cal F}_k(G)$ corresponds to a way of placing $k$ (indistinguishable) tokens in some of the $n$ vertices, say $1,2,\ldots,n$, of $G$.
For instance, $11122455\in CR_8^5$ corresponds to having 3 tokens in vertex 1, 2 tokens in vertex 2, no token in vertex 3, 1 token in vertex 4, and 2 tokens in vertex 5, which we represent by the sequence $32012$.
In general, we represent each vertex $\u$ of ${\cal F}_k(G)$ by a (non-negative) $n$-vector or $n$-sequence
$$
\u=(u_1,u_2,\ldots,u_n)\equiv u_1u_2\ldots u_n\quad \mbox{with } \sum_{i=1}^{n}u_i=k.
$$ 
 Moreover, two of its vertices, $\u$ and $\v$, are adjacent if one token of the multiset representing $\u$ is moved along an edge to another vertex (so obtaining the multiset representing $\v$). Then, notice that (the multisets of) $\u$ and $\v$ have $k-1$ elements in common.
Curiously enough, if $P_n$ is the path graph on $n$ vertices, it turns out that the $2$-supertoken graph ${\cal F}_2(P_n)$ is isomorphic to the 2-token graph $F_2(P_{n+1})$. The isomorphism between the vertices of ${\cal F}_2(P_n)$ and $F_2(P_{n+1})$ is given by the mapping $ij\mapsto i(j+1)$.

This kind of token graph was introduced by Hammack and Smith~\cite{HaSm17}, who named them \emph{reduced $k$-th power of graphs}, and provided
a construction of minimum cycle bases for them.

\begin{figure}[t]
	\begin{center}
		\includegraphics[width=14cm]{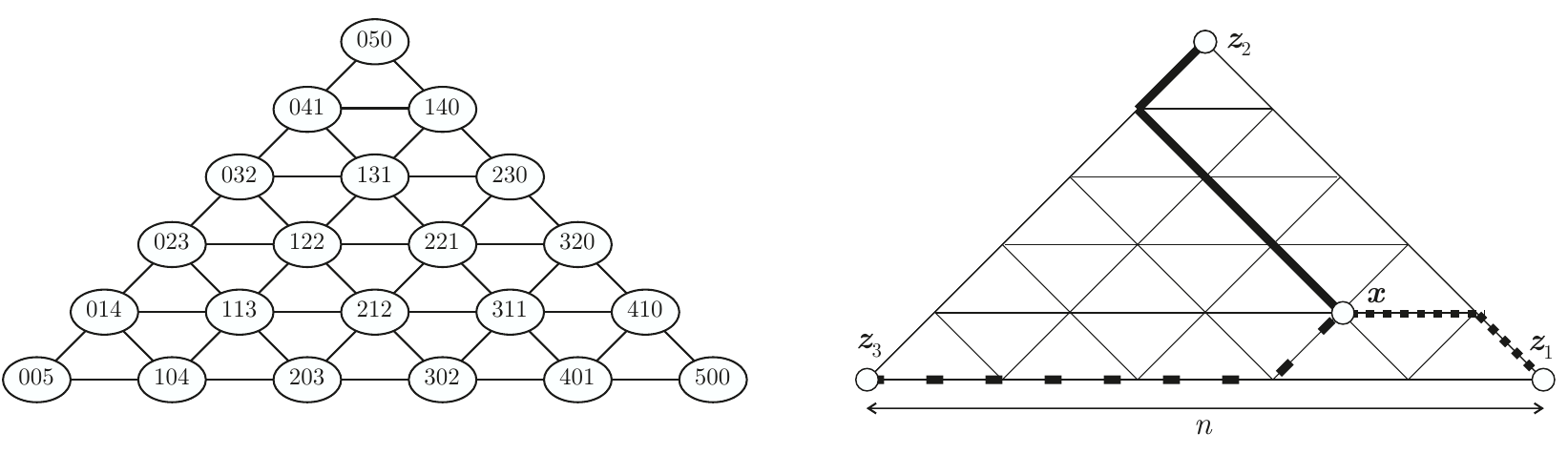}
	\end{center}
	\vskip-.4cm
	\caption{$(a)$ The supertoken graph ${\cal F}_5(K_3)$; $(b)$ In ${\cal F}_5(K_3)$ the shortest paths from vertex $\x$ to the vertices $\z_1=500,\z_2=050,\z_3=005$ are in short dashed line,  continuous line, and long dashed line, respectively. The position of $\x$ is given by the vector $\vecrho=(2,4,4)$.}
	\label{fig:ST_3,5_}
\end{figure}

\subsection{The case of ${\cal F}_k(K_n)$}

First, we consider the case of $G=K_n$, the complete graph on $n$ vertices.
As an example, ${\cal F}_5(K_3)$ is drawn in Figure \ref{fig:ST_3,5_}$(a)$.

In the following result, we describe some basic metric properties of ${\cal F}_k(K_n)$. 
\begin{Proposition}
\label{propo:Fk(Kn)}
Given integers $n,k$, the supertoken graphs ${\cal F}_k={\cal F}_k(K_n)$ satisfy the following statements. 
\begin{itemize}
\item[$(i)$]
The distance between two generic vertices $\x=x_1x_2\ldots x_n$ and $\y=y_1y_2\ldots$ $y_n$ is
\begin{equation}
\dist_{{\cal F}_k}(\x,\y)=\frac{1}{2}\sum_{i=1}^{n}|x_i-y_i|.
\label{dist(x,y)}
\end{equation}
\item[$(ii)$]
The eccentricity of a vertex $\x=x_1x_2\ldots x_n$ is
\begin{equation}
 \ecc_{{\cal F}_k}(\x)=k-\min_{1\le i\le n}{x_i}.
 \label{ecc(x)}
\end{equation}
\item[$(iii)$]
The diameter of ${\cal F}_k(K_n)$ is $\diam({\cal F}_k(K_n))=k$.
\item[$(iv)$]
The radius of ${\cal F}_k(K_n)$ is $\rad({\cal F}_k(K_n))=n-\lfloor \frac{n}{k}\rfloor$. 
\end{itemize}
\end{Proposition}
\begin{proof}
$(i)$
Let us consider the sets $X=\{i: x_i> y_i\}$, 
$Y=\{i: x_i< y_i\}$, and $Z=\{i: x_i=y_i\}$.
Then,
\begin{align}
\sum_{i\in X} (x_i-y_i) &=\sum_{i\in X} x_i-
\sum_{i\in X} y_i=n-\sum_{i\in Z} x_i-\sum_{i\in Y} x_i-\left(n-\sum_{i\in Z} y_i-\sum_{i\in Y} y_i\right) \nonumber\\
 &=\sum_{i\in Y} (y_i-x_i).\label{equality}
 \end{align}
 Then, to go from $\x$ to $\y$, in each step, we move a token from a vertex $i\in X$ (so reducing the value of $x_i$) to a vertex $j\in Y$ (increasing the value of $x_j$). Thus, the length of the path is $\sum_{i\in X}(x_i-y_i)$ and, from \eqref{equality}, $(i)$ follows.\\
 $(ii)$ Given the vertex $\x=x_1x_2\ldots x_n$ with $x_h=\min_{1\le i\le n} x_i$, the vertex $\y=y_1y_2\ldots y_n$ is at maximum distance from $\x$ when we have to move the maximum possible number of tokens, one at a time. That is,  when $y_j=0$ for $j\neq h$ and $y_h=k-x_h$.\\
 $(iii)$ According to \eqref{dist(x,y)}, the maximum possible distance between two vertices $\x$ and $\y$ is $k$. That is, when the $\sum_{i=1}^{n}|x_i-y_i|=2k$. Using the notation in the proof of $(i)$, this occurs when, for all $i\in X$, $\sum_{i\in X}x_i=k$ and $y_i=0$, (and, hence, for all $j\in Y$, $\sum_{j\in Y}y_j=k$ and $x_j=0$). \\
 $(iv)$ From \eqref{ecc(x)}, the minimum eccentricity of a vertex $\x$ is obtained when  $x_h=\min_{1\le i\le n}{x_i}$ is maximum. Since $\sum_{i=1}^n x_i=k$, this occurs when $x_h=\lfloor k/n\rfloor$.
    \end{proof}
    \begin{Example}
    In the graph ${\cal F}={\cal F}_3(K_5)$ of Figure \ref{fig:ST_3,5_}$(a)$, we get:\begin{itemize}
    \item[$(i)$]
    $\dist_{{\cal F}}(203,140)=\frac{1}{2}(1+4+3)=4$;
   \item[$(ii)$]
   $\ecc_{\cal F}(122)=5-\min\{1,2\}=4$;
  \item[$(iii)$]
  $\dist_{\cal F}(500,041)=\diam({\cal F})=5$;
   \item[$(iv)$]
   $\rad({\cal F}) =\ecc(122) = 5- \lfloor5/3\rfloor = 4$.
   \end{itemize}
   \end{Example}



\subsection{The general case of ${\cal F}_k(G)$}

Let $G=(V,E)$ be a (connected) graph with vertex set $\{1,2,\ldots,n\}$. The metric properties in the following result are direct generalizations of Proposition
\ref{propo:Fk(Kn)}.

\begin{Proposition}
Let $G$ be a graph on $n$ vertices, diameter $d$, and radius $r$. Given an integer $k$, the supertoken graphs ${\cal F}_k(G)$ satisfy the following statements. 
\begin{itemize}
\item[$(i)$]
The diameter of ${\cal F}_k(G)$  is $\diam({\cal F}_k(G))=kd$. 
\item[$(ii)$]
The radius of ${\cal F}_k(G)$ is $\rad({\cal F}_k(G))\le kr$. 
\item[$(iii)$]
If $G$ is an $r$-antipodal graph with diameter $d$, then ${\cal F}_k(G)$ has $r$ vertices that are mutually at distance $kd$. 
\end{itemize}
\end{Proposition}
\begin{proof}
    $(i)$ Let $i,j\in V(G)$ such that $\dist_G(i,j)=d$. Then, in ${\cal F}_k(G)$, to go from a vertex $\x$ to vertex $\y$ such that $x_i=k$ (with $x_h=0$ for $h\neq i$) and $y_j=k$  (with $y_h=0$ for $h\neq j$) corresponds to move, in $G$, $k$ tokens from vertex $i$ to vertex $j$, and this requires a maximum of $kd$ steps. \\
    $(ii)$ Let $i$ such that $\ecc(i)=r$. Then, we claim that, in ${\cal F}_k(G)$, the eccentricity of a vertex $\x$ such that $x_i=k$ (with $x_h=0$ for $h\neq i$) is at most $kr$. Indeed, in $G$, we can move every token in vertex $i$ to any vertex $j$ in, at most, $r$ steps. This means that we can go, in ${\cal F}_k(G)$, from $\x$ to any vertex $\y$ in, at most, $kr$ steps.\\
    $(iii)$ If $i_1,i_2,\ldots,i_r$ are vertices in $G$ at distance $d$ from each other, the $r$ vertices of ${\cal F}_k(G)$ $\x_1,\x_2,\ldots,\x_r$ such that $(\x_1)_{i_1}=k$, $(\x_2)_{i_2}=k$, $\ldots$,  $(\x_r)_{i_r}=k$  are, from $(i)$, mutually at distance $kd$.
    \end{proof}

Contrary to the case when $G=K_n$, when $G$ is a general graph, there is no closed formula for the distance between two vertices $\x$ and $\y$. Instead, we give the following  algorithm to find a shortest path from $\x$ to $\y$ in ${\cal F}_k(G)$:
\begin{enumerate}
    \item Given 
    $\x=x_1x_2\ldots x_n$ and $\y=y_1y_2\ldots y_n$, let $X=\{i:x_i>y_i\}$ and $Y=\{j:y_j>x_j\}$.
      Thus, if $x_h=y_h$ for some $h\in[1,n]$, such index $h$ is neither present in $X$ nor in $Y$.
    \item 
Consider the vectors 
\begin{itemize}
\item 
$\x'=(i_1^{x_{i_1}-y_{i_1}},i_2^{x_{i_2}-y_{i_2}},\ldots,i_{\sigma}^{x_{i_{\sigma}}-y_{i_{\sigma}}})$ such that $i_k\in X$ for $k=1,\ldots,\sigma$,
and $i_k^{x_{i_k}-y_{i_k}}$ stands for $i_k,\stackrel{(x_{i_k}-y_{i_k})}{\ldots\ldots}, i_k$,
\item
$\y'=(j_1^{y_{j_1}-x_{j_1}},j_2^{y_{j_2}-x_{j_2}},\ldots,j_{\sigma}^{y_{j_{\tau}}-x_{j_{\tau}}})$ such that $j_h\in Y$ for $h=1,\ldots,\tau$,
and $j_h^{y_{j_h}-x_{j_h}}$ stands for $j_h,\stackrel{(y_{j_k}-x_{j_h})}{\ldots\ldots}, j_h$.
\end{itemize}
(Notice that $\x'$ and $\y'$ can also be seen as multisets with the same number, say $n'=\sum_{i_k\in X}(x_{i_k}-y_{i_k})=\sum_{j_h\in Y}(y_{i_h}-x_{i_h})$, of vertices in $G$).
\item 
Define a distance $n'\times n'$ matrix $\D$ with rows and columns indexed by the entries of $\x'$ and $\y'$ respectively, and entries $(\D)_{i_k j_h}=\dist_G(i_k, j_h)$.
(This matrix represents a weighted complete bipartite graph $K^*_{n',n'}$, where the edge $\{i_k,j_h\}$ has weight 
$\dist_g(i_k,j_h)$).
  \item 
  Apply the known algorithm to find  a minimum weighted perfect matching  in $K^*_{n',n'}$ (see the Appendix).\\
  (This algorithm, usually called the `Hungarian method', is known to be `strongly polynomial'; that is, it is polynomial in the number of weights (distances) used. See, for instance, Edmonds \cite{e65} and Kuhn \cite{k55}.)
   \item 
   For each pair $i_k$ and $j_h$ of paired vertices of the matching, we move, in $G$, a token from
   $i_k$ to $j_h$.\\
   (Since in  $K^*_{n',n'}$ there are possibly vertices with the same labels, the above operation could be repeated several times).
\end{enumerate}

\begin{Example}
Let ${\cal{F}}_9(G)$ with $G=C_6$ is the cycle graph with vertex set $\{1,2,\ldots,6\}$. Let us consider the vertices $\x=310212$ and $\y=201132$ of ${\cal{F}}_9(G)$.
Then, since we only need to move the tokens that do not coincide between $\x$ and $\y$, $X=\{1,2,4\}$, $Y=\{3,5\}$, $\x'=(1,2,4)$, $\y'=(3,5^2)=(3,5,5)$. 
The distance matrix between the entries (vertices in $G$) of $\x'$ and $\y'$ is
\begin{equation*}
\D:=\begin{blockarray}{cccc}
& {\scriptstyle 3} & {\scriptstyle 5} & {\scriptstyle 5}  \\
\begin{block}{c(ccc)}
{\scriptstyle 1 } & 2 & 2 & \mathbf{2}  \\
{\scriptstyle 2} & \mathbf{1} & 3 & 3  \\
{\scriptstyle 4} & 1 & \mathbf{1} & 1  \\
\end{block}
\end{blockarray}\ ,
\qquad\qquad
\vcenter{\hbox{\includegraphics[width=2cm,height=2cm]{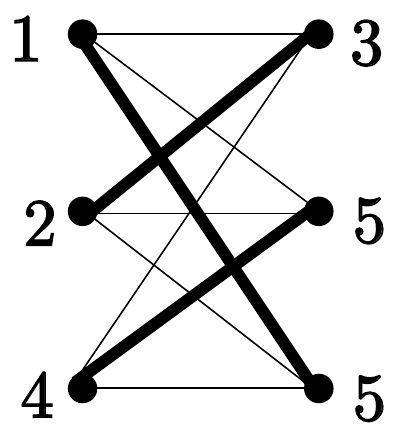}}}
\end{equation*}
where a perfect matching with minimum weight $4$ in $K_{3,3}$ consists of the edges  $\{2,3\}$,  $\{4,5\}$, and $\{1,5\}$ (see the above figure). This matching corresponds to the entries in boldface in $\D$.
Then, a path from $\x$ to $\y$ is (where `$i\rightarrow j$' means that $x_i\rightarrow x_i-1$
and $x_j\rightarrow x_j+1$ for $i\in X$ and $j\in Y$):
$$
310212\ \stackrel{2\rightarrow 3}{\sim}\ 301212\ \stackrel{4\rightarrow 5}{\sim}\ 301122\ \stackrel{1\rightarrow 5} {\sim}\ 201132,
$$
indicating that $\dist_{{\cal F}_9}(\x,\y)=1+1+2=4$, according to matrix $\D$.
\end{Example}

Notice that the minimization of the sum of the weights in $\D$ (or distances in $G$) ensures that this algorithm gives a shortest path.

For a specific resolving set of ${\cal F}_k(G)$, the position of a vertex $\x$ is easily computed.
Indeed, consider the distance between a generic vertex $\x=x_1x_2\ldots x_n$ and the vertex $\y=0\ldots0k0\ldots0$, where the $k$ is the $j$-th entry for some $1\le j\le n$. Then, 
$$
\x'=(1^{x_1},\ldots,(j-1)^{x_{j-1}}, (j+1)^{x_{j+1}},\ldots, n^{x_n})\quad\mbox{and}\quad \y'=(j^{k-x_{j}}).
$$
Then, the distance matrix $\D$ in point 3 of the algorithm has dimensions $(k-x_{j})\times (k-x_{j})$
 with columns indexed by $j,j,\stackrel{(k-x_j)}{\ldots},j$ and all its rows equal to
$$
(\dist(1,j)^{x_1},\ldots,\dist(j-1,j)^{x_{j-1}},\dist(j+1,j)^{x_{j+1}},\ldots, \dist(n,j)^{x_n}).
$$
Thus, we can pair $1$ with $j$ $x_1$ times, $2$ with $j$ $x_2$ times, etc. In other words, seen in $G$,
we move $x_1$ tokens from 1 to $j$ (with a total number of steps $x_1\dist(1,j)$), $x_2$ tokens from $2$ to $j$  (with a total number of steps $x_2\dist(2,j)$), etc.
Putting all together, to go from vertex $\x$ to $\y$ in ${\cal F}_k(G)$  requires $x_i\dist_G(i,j)$ steps for every $i=1,2,\ldots,n$. Thus, we have proved the following result.
\begin{Lemma}
\label{lem:dist(x,z)}
Given the supertoken graph ${\cal F}_k={\cal F}_k(G)$ of the graph $G=(V,E)$ with vertex set $V({\cal F}_k)=\{1,2,\ldots,n\}$, let us consider the subset of $V({\cal F}_k)$
\begin{equation}
    C=\{\z_1,\z_2,\ldots,\z_n\}=\{k0\ldots0,\ 0k0\ldots0,\ \ldots,\ 0\ldots0k\}.
    \label{resolv}
\end{equation}
Then, the distance between a vertex $\x=x_1x_2\ldots x_n\in V({\cal F}_k)$ and the vertex $\z_j$ is
\begin{equation}
\dist_{{\cal F}_k}(\x,\z_j)=\sum_{i=1}^n x_i\dist_G(i,j).
\end{equation}
\end{Lemma}
\begin{Corollary}
\label{coro:rho=xD}
  Let $G$ be a graph on $n$ vertices, with diameter $d$ and distance matrix $\D$.
  Let $C$ be the vertex subset of ${\cal F}_k={\cal F}_k(G)$ in Lemma \ref{lem:dist(x,z)}.
  Then, the position $\vecrho$ of a vertex  $\x=(x_1,x_2,\ldots,x_n)\in V({\cal F}_k)$ (represented as a vector) with respect to $C$ is
  \begin{equation}
       \vecrho=\vecrho(\x|C)=\x\D.
  \end{equation}
\end{Corollary}
\begin{proof}
Just notice that, for every $j=1,2,\ldots,n$,
    $$
  (\x\D)_j=\sum_{i=1}^nx_i(\D)_{ij}=\sum_{i=1}^nx_i\dist_G(i,j)=\dist_{{\cal F}_k}(\x,\z_j)=\rho_j,
    $$
    where we used the definition of $\D$ and Lemma \ref{lem:dist(x,z)}.
\end{proof}
The above results suggest the following definition.
\begin{Definition}
 Given a graph $G$ on $n$ vertices and an integer $k\ge 1$, we say that a nonnegative vector $\vecrho=(\rho_1,\rho_2,\ldots,\rho_n)$ is $(G,k)$-feasible if there exists a vertex $\x$ of ${\cal F}_k(G)$ whose position with respect to the set $C$ in \eqref{resolv} is $\vecrho$.   
\end{Definition}
In the case when the distance matrix $\D$ is nonsingular, we have the following characterization of feasible vectors.
\begin{Lemma}
Let $G$ have a nonsingular distance matrix $\D$. Then, a vector $\vecrho$ is $(G,k)$-feasible for some $k$ if and only if $\x=\vecrho\D^{-1}$ is a non-negative integer vector whose components sum up to $k$.
\end{Lemma}
\begin{proof}
    Clearly, the vector $\x$ so obtained represents a vertex of the supertoken graph ${\cal F}_k(G)$ if and only if it satisfies the conditions. Moreover, its position with respect to $C$ is $\x\D=\vecrho\D^{-1}\D=\vecrho$, as required.
\end{proof}
\begin{Example}
\label{D(Kn)}
The distance matrix of the complete graph $K_n$ coincides with its adjacency matrix $\D=\A=\Circ(0,1,1,\ldots,1)$, with inverse $\D^{-1}=\frac{1}{n-1}\Circ (2-n,1,1,\ldots,1)$.
Then, in the graph ${\cal F}_3(K_5)$ of Figure \ref{fig:ST_3,5_}$(a)$, we show that the vector $\vecrho=(2,4,4)$ is feasible since $\vecrho\D^{-1}=(3,1,1)$ corresponds indeed to a vertex of  ${\cal F}_3(K_5)$, see again Figure \ref{fig:ST_3,5_}$(a)$. In contrast, the vector $\rho=(1,3,3)$ is not feasible since $\vecrho\D^{-1}=\frac{1}{2}(5,1,1)$. 
\end{Example}

\section{On the metric dimension of ${\cal F}_k(G)$}
\label{sec:metric-dim}

\subsection{The case $G=K_n$}

We first present a result giving a bound on the metric dimension of ${\cal F}_k(K_n)$.

\begin{Proposition}
\label{propo1}
The metric dimension of ${\cal F}_k(K_n)$ satisfies $\dim({\cal F}_k(K_n))$ $\leq n-1$.
\end{Proposition}

\begin{proof}
Let us first prove that ${\cal F}_k(K_n)$ has the resolving set $C$ in \eqref{resolv}. Then, the distances between a generic vertex $\x=x_1x_2\ldots x_n$ and the vertices of $C$ are
$$
\mathrm{dist}(\x,\z_i)=\sum_{j\neq i} x_j=k-x_i,
$$
for $i=1,2,\ldots,n$. Thus, different vertices have different distance vectors to $C$, and so 
$\dim({\cal F}_k(K_n))\leq n$.
Moreover, the sum of the entries of such a vector is a constant. Indeed,
$$
\sum_{i=1}^{n}\mathrm{dist}(\x,\z_i)=\sum_{i=1}^{n}\sum_{\stackrel{j=1}{j\neq i}}^{n} x_j=\sum_{i=1}^{n}(k-x_i)=nk-\sum_{i=1}^{n}x_i=nk-k=(n-1)k.
$$
This means that any $n-1$ entries of a distance vector to $C$, say $\mathrm{dist}(\x,\z_i)$ for $i=1,\ldots,n-1$, determine the remaining one
$\mathrm{dist}(\x,\z_n)$. Hence, we conclude that $C'=\{\z_1,\z_2,\ldots,\z_{n-1}\}=\{k0\ldots0,$ $0k0\ldots0, 00k0\ldots0, 0\ldots k0\}$ is also a resolving set with $|C'|=n-1$, which proves the result.
\end{proof}


\begin{Conjecture}
The metric distance of the supertoken graph ${\cal F}_k(K_n)$ is
$$
\dim({\cal F}_k(K_n))= n-1.
$$
\end{Conjecture}

In support of this conjecture, we have the following result.

\begin{Proposition}
\label{propo2}
Given any fixed value of $n$, 
the metric distance of the supertoken graph ${\cal F}_k(K_n)$ is
$$
\dim({\cal F}_k(K_n))= n-1
$$
when $n$ is large enough.
\end{Proposition}

\begin{proof}
From Proposition \ref{propo1} and Lemma \ref{lemma1}, we only need to prove that we cannot have $\dim({\cal F}_k(K_n))\le n-2$, that is, that
\begin{equation}
\label{ineq1}
k^{n-2}+n-2<{n+k-1 \choose k}.
\end{equation}
Then, the result follows because, using Stirling's approximation $n!\sim\sqrt{2\pi n} (n/e)^n$, the binomial term turns to be, for large $n$,  of the order of $O(k^{n-1})$.
\end{proof}

For instance, for $n=3,4$, the inequation \eqref{ineq1} (and, hence, $\dim({\cal F}_k(K_n))= n-1$)  holds for all $k\ge 1$;  whereas for $n=5,6$, it holds for $k\not\in [5,10]$ and $k\not\in[3,104]$, respectively.

\subsection{The case of a general graph $G$}

Proposition \ref{propo1} is a particular case of the following result.

\begin{Theorem}
\label{th:dimFF}
Let $G=(V,E)$ be a graph with $n$ vertices, diameter $d$, and nonsingular distance matrix $\D$. Then, the metric dimension of the supertoken graph ${\cal F}_k(G)$ satisfies
\begin{equation}
\dim({\cal F}_k(G))\le |V|.
\label{dimFk(a)}
\end{equation}
Moreover, if $G$ is a degree-regular graph with sequence of degrees $k_1,k_2,\ldots,k_d$, 
then
\begin{equation}
\dim({\cal F}_k(G))\le |V|-1.
\label{dimFk(b)}
\end{equation}
\end{Theorem}
\begin{proof}
With $V=\{1,2,\ldots,n\}$, let us prove that the set in \eqref{resolv}
$$
C=\{\z_1,\z_2,\ldots,\z_{n}\} =\{k0\ldots0,\ 0k0\ldots0,\ \ldots,\ 00\ldots 0k\}
$$ 
is a resolving set.
With this aim, assume that two vertex labels of ${\cal F}_k(G)$, say $\x=(x_1,x_2,\ldots, x_n)$ and $\y=(y_1,y_2,\ldots, y_n)$, with $\sum_{i=1}^n x_i=\sum_{y=1}^n=k$ and $x_i,y_i\in\{0,\ldots,k\}$, yield the same distance vector. Then, by Corollary \ref{coro:rho=xD},
$$
\x\D=\y\D.
$$
Since $\D$ is nonsingular, this implies that $\x=\y$. Therefore, different vertices have different distance vectors and, hence, \eqref{dimFk(a)} follows.\\
Moreover, is $G$ is degree-regular, its distance matrix $\D$ has constant row sum $\lambda=\sum_{i=0}^d ik_i$. In other words, $\D$ has the eigenvalue $\lambda$ with the corresponding (left or right) eigenvector $\j=(1,1,\ldots,1)$. Then, the sum of the entries of every distance vector is a constant because
$$
\x\D\j^{\top}= \lambda\x\j^{\top}=\lambda\sum_{i=1}^n x_i=\lambda k.
$$
This means that any $n-1$ entries of a distance vector to $C$, say $\mathrm{dist}(\x,\z_i)$ for $i=1,\ldots,n-1$, determine the remaining one
$\mathrm{dist}(\x,\z_n)$. Hence, we conclude that 
$$
C'=\{\z_1,\z_2,\ldots,\z_{n-1}\}=\{k0\ldots0,\ 0k0\ldots0,\ \ldots, \ 0\ldots k0\}
$$ 
is also a resolving set with $|C'|=n-1$ vertices, which proves the result in \eqref{dimFk(b)}.
\end{proof}

We already showed some families of graphs $G$ satisfying this theorem, that is, with a nonsingular distance matrix.  Indeed, by Theorem  \ref{th:detD}, this holds if $G$ is a tree or a unicyclic graph whose unique cycle has odd order. Moreover, in Example \ref{D(Kn)}, we dealt with the case when $G$ is a complete graph. In the following example, we check out the theorem when $G=C_5$.

\begin{figure}[t]
	\begin{center}
		\includegraphics[width=6cm]{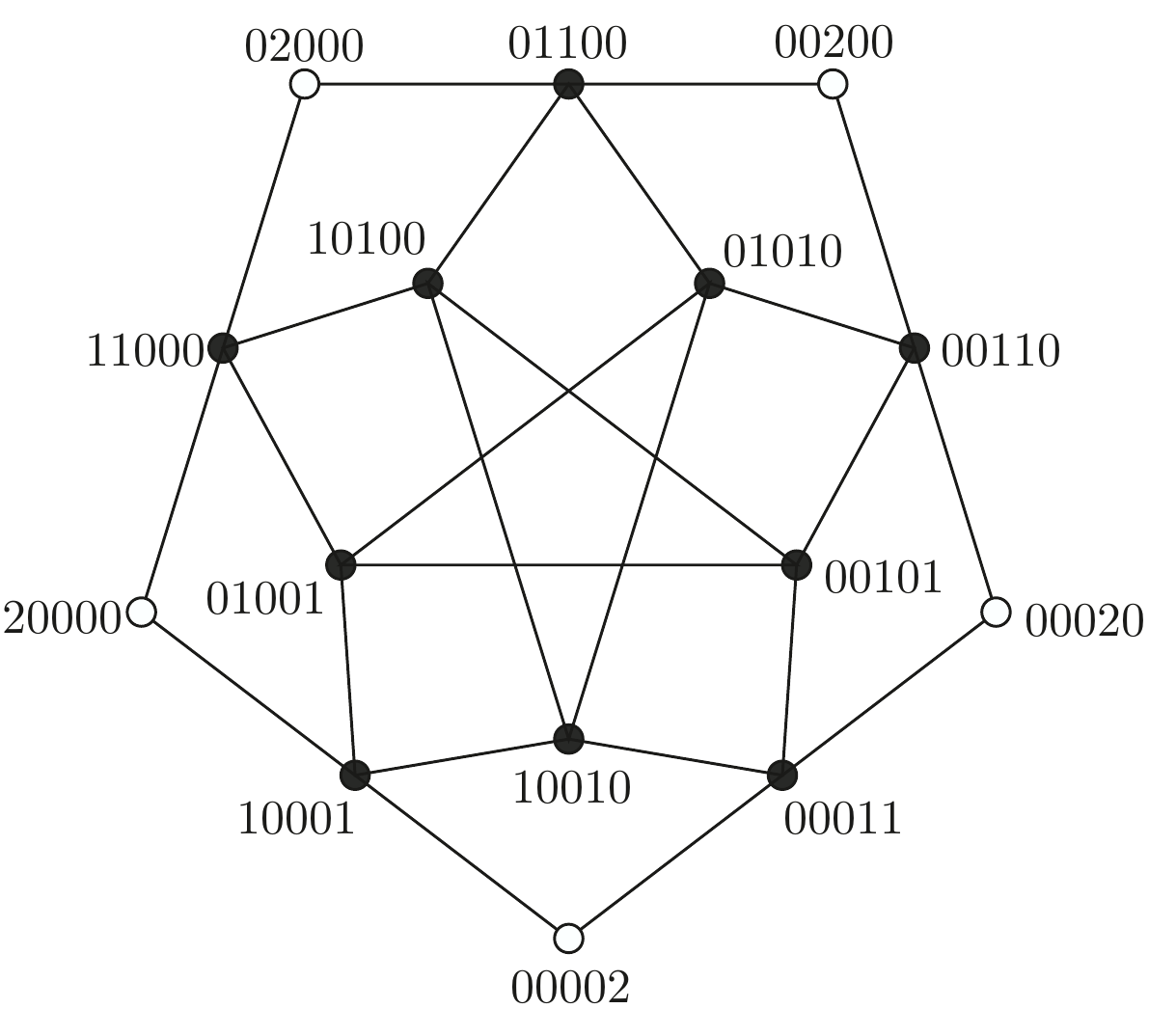}
	\end{center}
	\vskip-.25cm
	\caption{The supertoken graph ${\cal F}_2(C_5)$. The white vertices constitute a resolving set.}
	\label{fig:F^+_2(C5)}
\end{figure}

\begin{Example}
	Consider the graph $G=C_5$ and its supertoken ${\cal F}_2(C_5)$ (the latter is represented in Figure \ref{fig:F^+_2(C5)}). One can check that $\dim(C_5)=3$, whereas $\dim({\cal F}_2(C_5))=4$. The distance matrix of $C_5$ is the following circulant matrix:
 $$
 \D=\Circ(0,1,2,2,1)=\left(
 \begin{array}{ccccc}
     0 & 1 & 2 & 2 & 1\\
     1 & 0 & 1 & 2 & 2\\
     2 & 1 & 0 & 1 & 2\\
     2 & 2 & 1 & 0 & 1\\
     1 & 2 & 2 & 1 & 0
 \end{array}
 \right).
 $$
 Thus, by using Theorem \ref{th:detD}$(ii)$, we see that $\D(C_5)$ is nonsingular and, consequently, by Theorem \ref{th:dimFF}, $C=\{20000,02000,00200,00020,00002\}$ is a resolving set of ${\cal F}_2(C_5)$. Thus, let us check that the positions of vertices $\x=x_1x_2x_3x_4x_5=(x_1,x_2,x_3,x_4,x_5)$ with respect to $C$, that is, $\vecrho(\x)=\x\D$, where $x_i\in\{0,1,2\}$ and $\sum_{i=1}^5 x_i=2$, are all different. (The five vertices $\x$ with the same pattern, up to right shift, are the rows of the matrices on the left): 
 \begin{align*}
  & \left(
 \begin{array}{ccccc}
     2& 0 & 0 & 0 & 0\\
     0 & 2 & 0 & 0 & 0\\
     0 & 0 & 2 & 0 & 0\\
     0 & 0 & 0 & 2 & 0\\
     0 & 0 & 0 & 0 & 2
 \end{array}
 \right)\left(
 \begin{array}{ccccc}
     0 & 1 & 2 & 2 & 1\\
     1 & 0 & 1 & 2 & 2\\
     2 & 1 & 0 & 1 & 2\\
     2 & 2 & 1 & 0 & 1\\
     1 & 2 & 2 & 1 & 0
 \end{array}
 \right)=\left(
 \begin{array}{ccccc}
     0 & 2 & 4 & 4 & 2\\
     2 & 0 & 2 & 4 & 4\\
     4 & 2 & 0 & 2 & 4\\
     4 & 4 & 2 & 0 & 2\\
     2 & 4 & 4 & 2 & 0
 \end{array}
 \right), \\
 & \left(
 \begin{array}{ccccc}
     1 & 1 & 0 & 0 & 0\\
     0 & 1 & 1 & 0 & 0\\
     0 & 0 & 1 & 1 & 0\\
     0 & 0 & 0 & 1 & 1\\
     1 & 0 & 0& 0 & 1
 \end{array}
 \right)\left(
 \begin{array}{ccccc}
     0 & 1 & 2 & 2 & 1\\
     1 & 0 & 1 & 2 & 2\\
     2 & 1 & 0 & 1 & 2\\
     2 & 2 & 1 & 0 & 1\\
     1 & 2 & 2 & 1 & 0
 \end{array}
 \right)=\left(
 \begin{array}{ccccc}
     1 & 1 & 3 & 4 & 3\\
     3 & 1 & 1 & 3 & 4\\
     4 & 3 & 1 & 1 & 3\\
     3 & 4 & 3 & 1 & 1\\
     1 & 3 & 4 & 3 & 1
 \end{array}
 \right), \\
 & \left(
 \begin{array}{ccccc}
     1 & 0 & 1 & 0 & 0\\
     0 & 1 & 0 & 1 & 0\\
     0 & 0 & 1 & 0 & 1\\
     1 & 0 & 0 & 1 & 0\\
     0 & 1 & 0 & 0 & 1
 \end{array}
 \right)\left(
 \begin{array}{ccccc}
     0 & 1 & 2 & 2 & 1\\
     1 & 0 & 1 & 2 & 2\\
     2 & 1 & 0 & 1 & 2\\
     2 & 2 & 1 & 0 & 1\\
     1 & 2 & 2 & 1 & 0
 \end{array}
 \right)=\left(
 \begin{array}{ccccc}
     2 & 2 & 2 & 3 & 3\\
     3 & 2 & 2 & 2 & 3\\
     3 & 3 & 2 & 2 & 2\\
     2 & 3 & 3 & 2 & 2\\
     2 & 2 & 3 & 3 & 2
 \end{array}
 \right).
 \end{align*}
 As claimed, the positions of the vertices (rows of the matrices on the right) are all distinct.
\end{Example}

Going back to the proof of Theorem \ref{th:dimFF}, if the distance matrix $\D$ is singular, we cannot assure that the set $C$ in \eqref{resolv} is a resolving set. To illustrate this fact, we present the following example.

\begin{Example}
	Consider the graph $G=C_6$ and its supertoken ${\cal F}_2(C_6)$. The cycle $C_6$ has vertices $1,2,\ldots,6$ with  distance matrix
 $$
 \D=\Circ(0,1,2,3,2,1)=\left(
 \begin{array}{cccccc}
     0 & 1 & 2 & 3 & 2 & 1\\
     1 & 0 & 1 & 2 & 3 & 2\\
     2 & 1 & 0 & 1 & 2 & 3\\
     3 & 2 & 1 & 0 & 1 & 2\\
     2 & 3 & 2 & 1 & 0 & 1\\
     1 & 2 & 3 & 2 & 1 & 0
 \end{array}
 \right).
 $$
 Applying again Theorem \ref{th:detD}$(ii)$, we observe that $\det\D=0$ and, hence, $\D$ is singular.
 Thus, if we take the  set $C=\{200000,020000,\ldots,000002\}$, the position of a vertices $\x_1=14=(1,0,0,1,0,0)$, $\x_2=25=(0,1,0,0,1,0$, and $x_3=36=(0,0,1,0,0,1)$ with respect to $C$ turns out to be the same
 $$
 \vecrho(x_i|C)=\x_i\D=(3,3,3,3,3,3)\quad \mbox{for $i=1,2,3$,}
 $$
 showing that $C$ is not a resolving set.
\end{Example}



\section{Appendix}
\label{sec:appendix}

 Let $G$ be a complete bipartite graph with independent sets $V_1$ and $V_2$, where the edge $\{i,j\}$ carries the weight $w_{ij}\ge 0$. Every matching $M$ of $G$ is represented by the incidence vector $\x$ with components $x_{ij}=1$ if $\{i,j\}\in M$, and $x_{ij}=0$ otherwise. Then, a solution of the following linear programming problem (LP) provides a minimum weighted perfect matching in $G$:
\begin{equation}
\boxed{
\begin{array}{rl}
 {\tt minimize} & \sum_{i\in V_1,j\in V_2} w_{ij}x_{ij}\\
 {\tt subject\ to} & \sum_{j\in V_2}x_{ij}=1,\quad i\in V_1\\
                       & \sum_{i\in V_1}x_{ij}=1,\quad j\in V_2\\
                       & x_{ij}\ge 0,\qquad i\in V_1, j\in V_2\\
\end{array}
}
\label{LP}
\end{equation}
Here, it is worth noting that the solutions to this problem are integers, $x_{ij}\in \{0,1\}$, and correspond to the vertices of a so-called `bipartite perfect matching polytope'. 
In fact, the algorithm also works for any bipartite graph since, by assigning infinite costs to the edges not present, one can assume that the bipartite graph is complete.

\end{document}